\documentclass[11pt]{article}

\usepackage{epsf}
\usepackage{psfrag}
\usepackage{amsmath}
\usepackage{amssymb}
\usepackage{lscape}
\usepackage{latexsym}
\usepackage[usenames]{color}
\usepackage[mathcal]{eucal}
\usepackage{stmaryrd}
\usepackage{textcomp}
\usepackage{setspace} 
\usepackage{graphicx} 
\usepackage{pgfplots} 
\usepackage{epsfig}   
\usepackage{epstopdf} 
\usepackage{appendix}
\usepackage{afterpage}
\usepackage{rotating} 
\usepackage{float} 
\usepackage[misc,geometry]{ifsym} 
\usepackage{comment}

\newtheorem{example}{Example}[section]

\newcommand{\ds}{\displaystyle}
\newcommand{\dZ}{{\cal Z \kern -0.7em Z}}
\newcommand{\dC}{{\rm\hbox{C \kern-0.8em\raise0.2ex\hbox{\vrule height5.4pt width0.7pt}}}}
\newcommand{\dQ}{{\rm\hbox{Q \kern-0.85em\raise0.25ex\hbox{\vrule height5.4pt width0.7pt}}}}

\newcommand{\proofbox}{\hspace{\fill}{$\Box$}}

\newcommand\old[1]{}
\newcommand{\beqa}{\begin{eqnarray*}}
\newcommand{\eeqa}{\end{eqnarray*}}

\begin{document}

\title{\bf \Large{Optimizing Inventory Management through Multiobjective Reverse Logistics with Environmental Impact}}
\author{I. B. Wadhawan \thanks{Centre for Smart Analytics (CSA), Institute of Innovation, Science and Sustainability, Federation University Australia, 
{\tt email: i.wadhawan@federation.edu.au}}
\and
M. M. Rizvi \thanks{Centre for Smart Analytics (CSA), Institute of Innovation, Science and Sustainability, Federation University Australia, {\tt email:  m.rizvi@federation.edu.au} (corresponding author).}}
\maketitle
\textbf{Abstract} We present novel mathematical models for inventory management within a reverse logistics system. Technological advancements, sustainability initiatives, and evolving customer behaviours have significantly increased the demand for repaired products. Our models account for varying demand levels for newly produced and repaired items. To optimize overall costs with constrained scenarios, we formulated mixed integer programming problems. Solution procedures for the proposed problems are introduced, and the accuracy of these solutions has been validated through numerical experiments. Additionally, we address the cost of waste disposal as an environmental concern. This paper develops a multiobjective mathematical model and provides an algorithm for the Pareto solution. Various scalarization techniques are utilized to identify the Pareto front, and a comparison of these techniques is presented.


\textbf{AMS} subject classifications. 90B05, 90C25, 90C29, 90C30.


\section{Introduction}\label{sec1}
Protecting the environment by recycling material and remanufacturing used products are inevitable options to reduce waste generation. Reducing waste can be attained through product recovery through refurbishing, repair, recycling and remanufacturing. Reverse logistics (RL) has been receiving increased attention because of the concrete benefits that flow from the repair and are seen as a solution to environmental waste problems. The research on coordinating repair with manufacturing and waste management was started in the 1970s but only limited research has considered the waste as a disposal cost \cite{Bazan2015}. 
\par
The earliest reported work is by Schrady \cite{Schrady1967}, who establishes a deterministic Economic Order Quantity (EOQ) model for repaired items with instantaneous manufacturing and recovery rates. Nahmias and Rivera \cite{Nahmias1979} extended Schardy's model for finite repair rate and one production cycle. Richter \cite{Richter1996a, Richter1996b} proposed an EOQ model and determined the optimal number of production and remanufacturing batch sizes based on the return rate. The system investigated by Richter \cite{Richter1997} consists of the serviceable and the repairable stock, which was governed by ``dispose all'' or ``recover all'' strategies. Dobos and Richter \cite{Dobos2004} generalised the work by Richter \cite{Richter1997} by considering finite production and repair rates. They extended the model by considering the quality of used products may not be suitable for recycling \cite{Dobos2006}. Teunter \cite{Teunter2001} conducted an analysis analogous to that of Richter \cite{Richter1997} comparing two strategies: (1) the maximum remanufacturing strategy, involving multiple recovery cycles and a single production (manufacturing) cycle per interval, and (2) the maximum production strategy, which incorporates multiple production cycles and a single recovery cycle per interval. These strategies were extended with scenarios featuring indefinite recovery times. El Saadany and Jaber \cite{Saadany2010} considered a mixed strategy of production and remanufacturing as opposite to Richter \cite{Richter1997}. Jaber and Saadany \cite{Jaber2009} proposed a model with the emphasis that remanufacturing items are considered as lower quality by the customers, which was extended by Hasanov et al \cite{Hasanov2012} by incorporating different cases of backorder. Saadany \cite{Saadany2013} has investigated remanufacturability issues of reverse logistics and changed a common assumption by differentiating unlimited remanufacturing of repairable products from the production of manufactured products. In addition, disposal cost is considered in the total cost function. Bazan \cite{Bazan2015} also considered waste cost as a part of the overall cost to minimize. Singh et al \cite{Singh2016} extended the work by Hasanov et al introducing finite production and remanufacturing rates. Bazan et al \cite{Bazan2016} reviewed the mathematical inventory model and established a production reliability model with random demand. Kozlovskaya et al \cite{Kozlovskaya2017} consider a general class of deterministic inventory models with constant demand and return to minimise the costs involved (sum of the production cost, holding cost and switching cost). Singh et al \cite{Singh2020} considered a reverse logistic model with variable production and remanufacturing under learning effects. 
\par
Currently, government pressure and legislation have contributed to the increased motivation for sustainability practices, but limited research has been conducted on this. Forkan et al \cite{Forkan2022} further extended the work by Nahmias and Rivera {\cite{Nahmias1979} by considering two different demands for manufacturing and repaired items along with the environmental issues such as GHG and energy consumption in their model. 
\par
Our paper considers a mathematical approach to optimize the cost associated with reverse logistics with different demands and the cost associated with waste, hence converting the model to a multiobjective problem. This paper focuses on obtaining the optimal batch cycles by solving a Mixed Integer Programming (MIP) problem, which has not been considered in any research so far. The normal practice is fixing the number of cycles and optimizing the cost. However, our paper is obtaining optimum cycles along with the cost by solving the mixed integer programming problem. 
\par
In summary, the objectives of the current research are as follows.
\begin{itemize}
\item [(i)] Develop an efficient material flow model in the reverse logistics system to calculate average inventory cost where repaired items are not considered the same as new items.
\item [(ii)] Develop mathematical expressions to calculate the total cost functions (involving optimised lot sizes and number of cycles) of the reverse logistics system for the case when the demand for the repaired and new procurement items are not the same.
\item [(iii)] Minimise environmental impact by considering waste disposal while considering the reverse logistics system of inventory management.
\item [(iv)] Develop a model for the multiobjective functions with inventory cost and waste cost and then minimise these conflicting functions simultaneously by using different scalarization techniques for the proposed models. 
\item [(v)] Design efficient algorithms and establish a solution procedure for the proposed models and implement them through computational experiments.
\end{itemize}

The rest of this paper comprises the following sections. Section \ref{Form} describes some preliminaries and assumptions to develop the proposed model. Section \ref{Moo} develops formulating three new mathematical reverse logistic models for constrained optimization problems. A multiobjective model with a solution approach is proposed in this section. Section \ref{NuExp} comprises numerical experiments and sensitivity analysis. Section \ref{Con} provided the conclusion and further research directions. An alternative diagram of inventory behaviour for new, collected, used and repaired items to the proposed problem is presented in the Appendix. 
 
\newpage
\section{Formulations and Solutions of the Proposed Models }\label{Form}
The following parameters are used in the formulation of the model.
\begin{table} [H]
\begin{tabular}{l l}
$T$& Length of new procurement and repair cycle.\\
$T_p$& Length of new items cycle.\\
$T_r$& Length of repaired items cycle.\\
$Q_p$& Batch size for new procurement over the time $T_p/n$.\\
$Q_r$& Batch size for repaired items over the time $T_r/m$.\\
$D_p$& Demand rate for new items/primary market.\\
$D_r$& Demand rate for repaired items/secondary market.\\
$p$  & Proportion of $D_p$ of new items (units per unit time) $0<p\leq 1$.\\
$1-p$& Lost proportion of $D_p$ of new items (units per unit time).\\
$q$  & Collection proportion of available returns of $pD_p$ ($0\leq q \leq 1$). \\
$1-q$& Waste proportion of $pD_p$.\\
$r$  & Proportion of Dr of repaired items (units per unit time) $0\leq r\leq 1$.\\
$1-r$& Lost proportion of Dr of repaired items (units per unit time).\\
$s$  & Collection proportion of available returns of $rD_r$ ($0\leq s \leq 1$). \\
$1-s$& Waste proportion of $rD_r$.\\
$qpD_p+srD_r$& Repairable stock rate ($R_1 = qpD_p, R_2 = srD_r$).\\
$S_p$ and $S_r$ & Setup cost of supply and repair depot per unit time $T$.\\
$h_p$ and $h_r$& Holding cost of supply and repair depot respectively per unit time $T$.\\
$cw$& Waste cost per unit time $T$.
\end{tabular}
\end{table}

\par
The following assumptions are used to ensure the ease of the developed model.
\begin{itemize}
\item Repaired items are of a different quality than the new items therefore, the demand for the new and repaired items is not fixed in time and may not be the same \cite{Jaber2009, Singh2016}.
\item Shortage is not allowed at the stock point and is represented by the constraint $qpD_p+srD_r =R_1+R_2 \geq D_r$. 
\item Lead time is zero (instantaneous replenishment).
\item Returned items undergo an inspection before being accepted at the repair depot. 
\item We assume the procurement period is never smaller than the repair period $T_r \leq T_p$ and $T = T_p+T_r$.
\end{itemize}

\par
We consider two types of depots (primary/supply and secondary/repair) in our proposed model as depicted in Figure~\ref{figFlowModel}. The first depot is the supply depot, which stores the new procurement and repaired items. The supply depot can satisfy customers' demand for new $D_p$ and repaired $D_r$ items. When the items are repaired, first they are stored at the repair depot and then shipped to the supply depot. The portion of demand $(1-p)D_p$ and $(1-r)D_r$ never return from the primary and secondary market to the system, so we assume that it does not have any influence on the system. Therefore, $pD_p$ and $rD_r$ are the portions of the demand which is returned to the system for either repair or disposal. Returned items vary in quality, and items with quality below the acceptance threshold are rejected. In the process, the amount of the accepted items is $R_1+R_2 =qpD_p+srD_r$ and the failed/waste item is $(1-q)pD_p$ and $(1-s)rD_r$. This failed quantity is considered in our model as a second objective function to account for waste cost.   
\begin{figure}[h]
	\centering
	\includegraphics[width=1.00\textwidth]{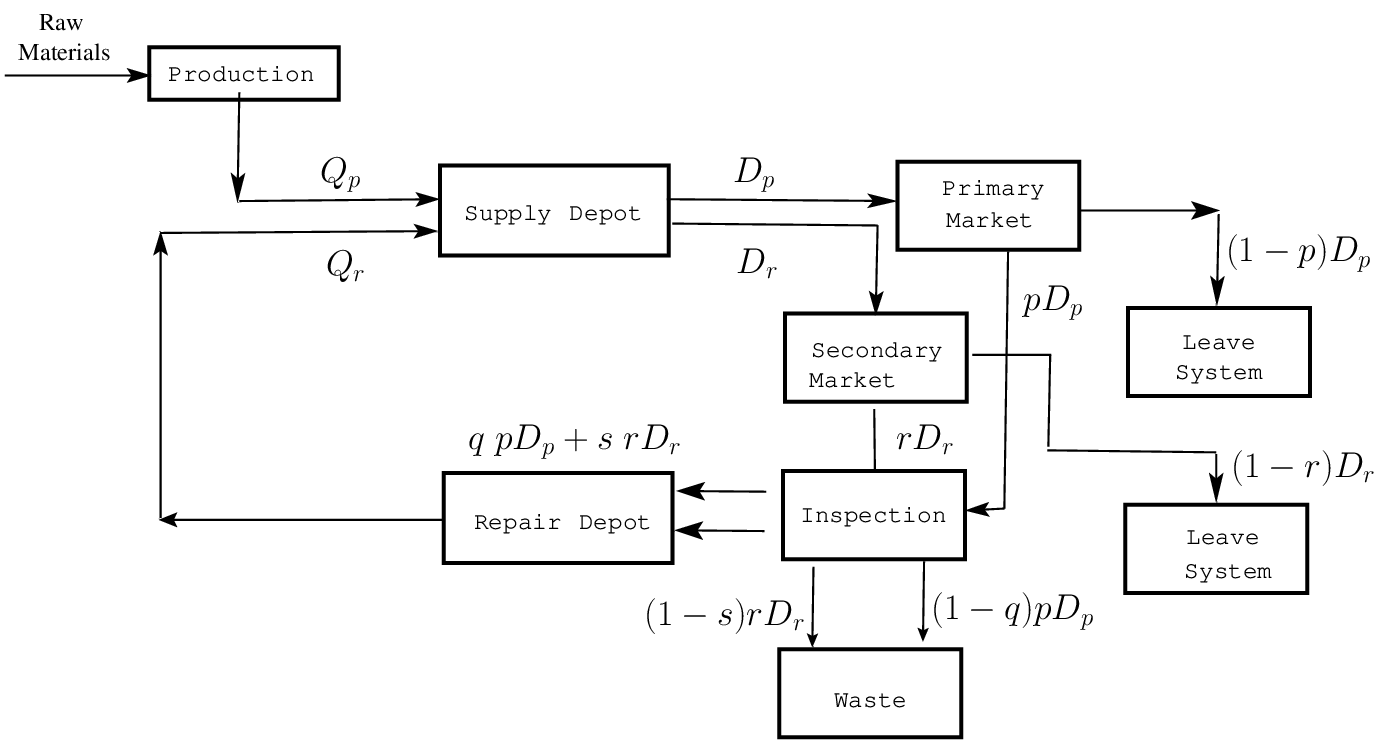}
	\caption{Material flow for a new procurement and repair system.}
	\label{figFlowModel}
\end{figure}

\par  
The inventory behaviour for new and repaired items is illustrated in Figure~\ref{figinventoryModel1}. The model begins with fulfilling the demand for newly produced and remanufactured items. As shown in Figure~\ref{figinventoryModel1}, the inventory of the supply depot always decreases at rate $D_p$ for new items and at rate $D_r$ at the repair depot for repaired items. In the proposed models, we consider $n$ cycles for new procurement and $m$ cycles for repair over the time $T = T_p+T_r$. 
\par
\noindent At the time $T=0$, since there are no items in the repair depot, only new procurement $Q_p$ is sold at the rate $D_P$. As new items are sold in $n$ cycles over a time period $T_p$, the inventory in the repair depot grows at the rate $pqD_p$. Therefore, we obtain $T_p D_p = n Q_p$. The same cycle repeats for the repaired items. The inventory of the repair depot increases at the rate $srD_r$ as we sell the repaired procurement $Q_r$ with the rate $D_r$ in $m$ cycles over period $T_r$. Therefore, we obtain $T_r D_r = m Q_r$.
\begin{figure}[ht!]
	\centering
	\includegraphics[width=1.00\textwidth]{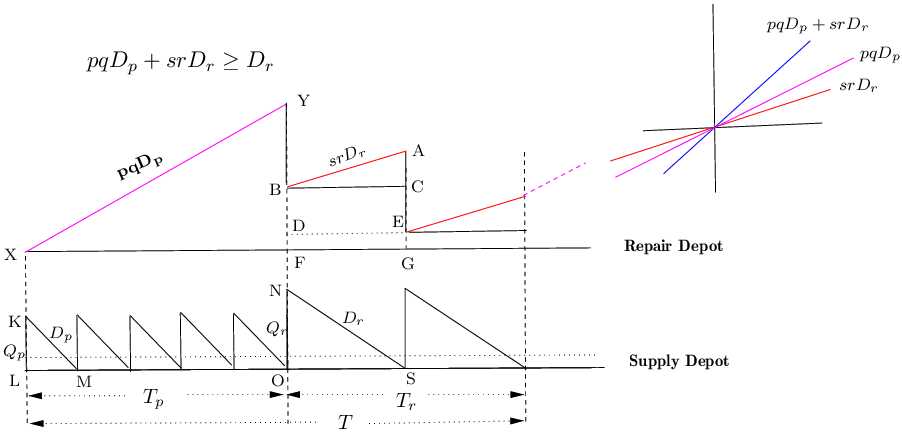}
	\caption{The behaviour of inventory for new, used and repaired items over the interval $T$.}
	\label{figinventoryModel1}
\end{figure}
The total cost per unit time for the proposed reverse logistic system $C(Q_p, Q_r, m,n)$ is the sum of the setup cost per unit time for the system, the total inventory holding cost per unit time in the supply depot and the repair depot. This is achieved by adding the areas in the supply depot ($n~\mbox{KLM}+m~\mbox{NOS}$) and the repair depot ($\mbox{XYF}+m~\mbox{ABC}+(m-1)~\mbox{BCED}+m~\mbox{DEGF}$) and is given as
\begin{multline}\label{eqcost1a}
C(Q_p, Q_r, m, n) = \frac{1}{T}\Bigg[ mS_n+nS_p + \left(\frac{Q_pT_p}{2}+\frac{Q_rT_r}{2}\right)h_p +\Bigg[ \frac{R_1T_p^2}{2} + \frac{R_2T_rQ_r}{2D_r} \\ +\frac{(m-1)}{2}T_r Q_r\left(1-\frac{R_2}{D_r}\right) 
+T_r\left(T_pR_1-Q_r-(m-1)\left(1-\frac{R_2}{D_r}\right)Q_r\right)\Bigg] h_r\Bigg],
\end{multline}
where, $T=T_p+T_r$, $R_1 = pqD_p$ and $R_2 = srD_r$.\\
For detailed information on area calculation, please check the Appendix section \ref{appen}.

\newpage
\section{Mathematical Model and Solution Approach}\label{Moo}
Numerous studies have addressed the disposal of non-repairable items \cite{Richter1996a, Richter1996b, Richter1997, Saadany2013}, with only a few considering it as a waste cost in the overall cost function \cite{Bazan2015}. However, none have examined this cost from an environmental perspective. Controlling this cost is crucial when optimizing the holding cost in reverse logistics systems. This paper explores the trade-offs between overall cost and waste cost, as well as between holding costs and waste cost.
\par 
\noindent The disposal cost function can be written as
\begin{equation}\label{eqcost1b}
f_2(s) = C_w\Big((1-q)pD_p+(1-s)rD_r\Big).
\end{equation}
Thus, we construct the reverse logistics problem as a mixed integer programming problem and propose three mathematical models. The objective of the first model is to minimize the inventory cost with the condition that the length of the new items cycle is always greater than the length of the repaired items cycle. In this particular model, the percentage of items $s$ approved for repair from the returned items from the secondary market remains constant. The new and repaired items are represented by integers, each requiring a minimum of one cycle.

\noindent {\bf Model 1: When $s$ is fixed}
\begin{equation}\label{eqmathemodel1}
\begin{array}{rl} 
\min & \ C(Q_p,Q_r,m,n)\\
\mbox{subject to the constraints} \\
& \ds T_r \leq T_p, \\ 
&\ds m, n \in {\mathbb{Z} \geq 1}.
\end{array}
\end{equation}
The second model also aims to minimize inventory costs with additional conditions. In addition to the first model, this model ensures that stock is always available and considers the changing percentages of items approved for repair from the returned items in the secondary market. We define \eqref{eqcost1a} as $f_1(Q_p,Q_r,m,n,s)$ when $S$ is considered as a variable.
\noindent {\bf Model 2: When $s$ is variable}
\begin{equation}\label{eqmathemodel2}
\begin{array}{rl} 
\min & \ f_1(Q_p,Q_r,m,n,s)\\
\mbox{subject to the constraints} \\
& \ds T_r \leq T_p, \\ 
& \ds D_r \leq R_1+R_2, \\
&\ds  m, n \in {\mathbb{Z} \geq 1},\\
& \ds 0 \leq s \leq 1.
\end{array}
\end{equation}
The third model considers two functions, inventory cost and waste cost while taking into account all the constraints from Model 2. The objective is to generate a set of solutions that minimize both the inventory cost and the environmental impact of waste cost.

\noindent {\bf Model 3: Multiobjective - Inventory cost versus Waste Cost}
\begin{equation}\label{eqmathemodel3}
\begin{array}{rl} 
\min & \ds \ \left\{f_1(Q_p,Q_r,m,n,s), f_2(s)\right\}\\
\mbox{subject to the constraints} \\
& \ds T_r \leq T_p, \\ 
& \ds D_r \leq R_1+R_2, \\
&\ds  m, n \in {\mathbb{Z} \geq 1},\\
& \ds 0 \leq s \leq 1.
\end{array}
\end{equation}

\subsection{Solution procedure}\label{solprocedure}
 \begin{description}
\item[Step $\mathbf{1}$] \textbf{(Input)} \\ Input the parameters $D_p$, $D_r$, $p$, $q$, $r$, $h_r$, $h_p$, $S_p$ and $S_r$. Fix $s$ for Model 1 with constraints represented by (\ref{eqmathemodel1}) as described in Example~\ref{ex1}, and variable $s$ for Model 2 with constraints presented in \eqref{eqmathemodel2} discussed in Examples~\ref{ex2} and \ref{ex3}.
 \item[Step  $\mathbf{2}$] \textbf{(Define Problem)}\\
 Define the cost function and constraints. Set the problem as a mixed integer problem (MIP) and execute MIP by using appropriate solvers such as Knitro, Bonmin, SCIP \cite{Achterberg2009} or BARON. For Model 1, the variables $Q_p$, $Q_r$, $m$ and $n$ are all defined as integers. However, for Model 2, $s$ is defined as a continuous variable. 
 \item[Step  $\mathbf{3}$] \textbf{(Solve Problem)}\\
Solve Problem $\ds \min f_1$, subject to the constraints of Model 1 represented by (\ref{eqmathemodel1}). The solver provides the solution $Q_p$, $Q_r$, $m$ and $n$ with optimized value of cost. 
 \item[Step  $\mathbf{4}$] \textbf{(Solver Validation)}\\
 Utilise different solvers such as Knitro, Bonmin, SCIP   and BARON in MATLAB environment \cite{Jonathan2014} and AMPL to solve MIP and to cross-check the optimum solution obtained. 
\end{description}
Repeat the same process with Model 2 with constraints represented by (\ref{eqmathemodel2}).

\section{Numerical experiments, results and sensitivity analysis}\label{NuExp}
This section provides numerical examples to demonstrate the reverse logistics models with various constraints. We used solution procedure \eqref{solprocedure} to optimize Models~$1$ and $2$ by using solvers such as Bonmin, BARON, and SCIP. As our proposed model is an MIP, We require solvers capable of efficiently handling such complexities that arise in MIP. We tested the above-listed solvers in Matlab and AMPL, and Bonmin appears to have been the most efficient in terms of computational time and accuracy in approximating the solution. Our analysis performs the computations on an HP Evo laptop equipped with 16 GB of RAM and a Core i7 processor running at 4.6 GHz.
\par
\noindent{\bf Case I (Model 1):} When $s$ is fixed 
\begin{example}\label{ex1} 
 Table 1 presents the relationship among costs, new item cycle $n$, repair cycle $m$ and time cycle $T$ when $D_p=100$, $D_r=43$, $s=0.7$, $q=0.8$, $h_r = 1.2$, and $S_r = 1$ respectively with $p$ and $r$ takes values from $40\%$ to $80\%$ and holding cost $h_p$ and setup cost $S_p$ switching between $1.6$ and $10$. These parameter values are taken from the Forkan et al \cite{Forkan2022} that belongs to the tyre industry. We aim to minimize Model 1 represented by (\ref{eqmathemodel1}). 
 \end{example}
	\begin{table} [H]
 		\caption{\small{\textit{Example-\ref{ex1} -- Numerical performance of model \eqref{eqcost1a} when $R_1+R_2\geq D_r$. }}}
		\footnotesize
		\vskip 1.5em
 		\centering
 		\begin{tabular}{|c| c| c| c| c| c| c| c| c| c| c| c|}
 			\hline
 			$p$ &   $r$ &  $h_p$  & $S_p$ & $Q_p$& $Q_r$ & $m$ &$n$ &{\small $C(Q_p,Q_r, m, n)$} & $T_r$& $T_p$ &  CPU time[sec] \\ [0.5ex]
 			\hline
0.4&0.4 &1.6 & 10 &35&15&3&3&51.78&1.05&1.05&6.2\\ [.5ex]\hline
0.4&0.4 &10 & 1.6 &5&3&15&21&62.81&1.05&1.05&3.3\\ [.5ex]\hline
0.6&0.6 &1.6 & 10 &35&15&3&3&68.78&1.05&1.05&4.2\\ [.5ex]\hline
0.6&0.6 &10 & 1.6  &5&3&15&21&79.81&1.05&1.05&1.05\\ [.5ex]\hline
0.8&0.8 &1.6 & 10 &35&15&3&3&85.78&1.05&1.05&4.2\\ [.5ex]\hline
0.8&0.8 &10 & 1.6  &5&3&15&21&96.81&1.05&1.05&4.11\\ [.5ex]\hline
                \end{tabular}
		\label{table:ex1}
	\end{table}
 \par
\noindent According to the assumptions of Model 1 in our experiment, the repairable stock rate $R_1+R_2$ is always more than the demand rate $D_r$. As observed in Table 1, when the repair stock rate increases while demand remains constant—meaning more items are received at the repair depot than are demanded—the optimal cost also rises. This is illustrated in the first, third, and fifth rows of the table, where the cost increases from $51.78$ to $85.78$ and, similarly, from $62.81$ to $96.81$ in rows 2, 4, and 6. Additionally, when the setup cost exceeds the holding cost, the overall cost is minimised, as seen in rows 1, 3 and 5, compared to rows 2, 4 and 6. 

\begin{example} \label{ex2}
Parameter values are taken from the Saadany \cite{Saadany2009}\\
$p=0.667$, $q=1$, $r=0.667$, $Dp=200$, $Dr=200$, $Sp=144$, $Sr=72$, $hp=12$ and $hr=3$.\\
CPU time = 46 sec
\end{example}

\begin{figure}[H]
\hspace{-1cm}
\begin{minipage}{130mm}
\hspace*{-2.5cm}
\begin{table}[H]
 		\caption{\small{\textit{Example- Numerical performance of Model 1}}} 
		\footnotesize
		\vskip 1em
 		\centering
 		\begin{tabular}{|c| c| c| c| c| c| c| c|}
 			\hline
 			 $s$ & $Q_p$& $Q_r$ & $m$ &$n$ &{\tiny $C(Q_p,Q_r, m, n)$} & $T_p$& $T_r$  \\ [0.5ex]
 			\hline
 		0.5  &	67	&50 & 4 & 3 & 829.6 &1.005&1\\ [.5ex]\hline
            0.7  &	67	&50 & 4 & 3 & 894.5 &1.005&1\\ [.5ex]\hline
            0.99 & 67	&50 & 4 & 3 & 923.4 &1.005&1\\ [.5ex]\hline
            0.995 & 67	&50 & 4 & 3 & {\bf 923.9} &1.005&1\\ [.5ex]\hline
            1 & 67&50 & 4 & 3 & 924.4 &1.005&1\\ [.5ex]\hline
		\end{tabular}
		\label{table:exp-1a}
	\end{table}
\end{minipage}
\hspace*{-2.5cm}
\begin{minipage}{50mm}
\vspace*{1.4cm} 
\includegraphics[width=50mm]{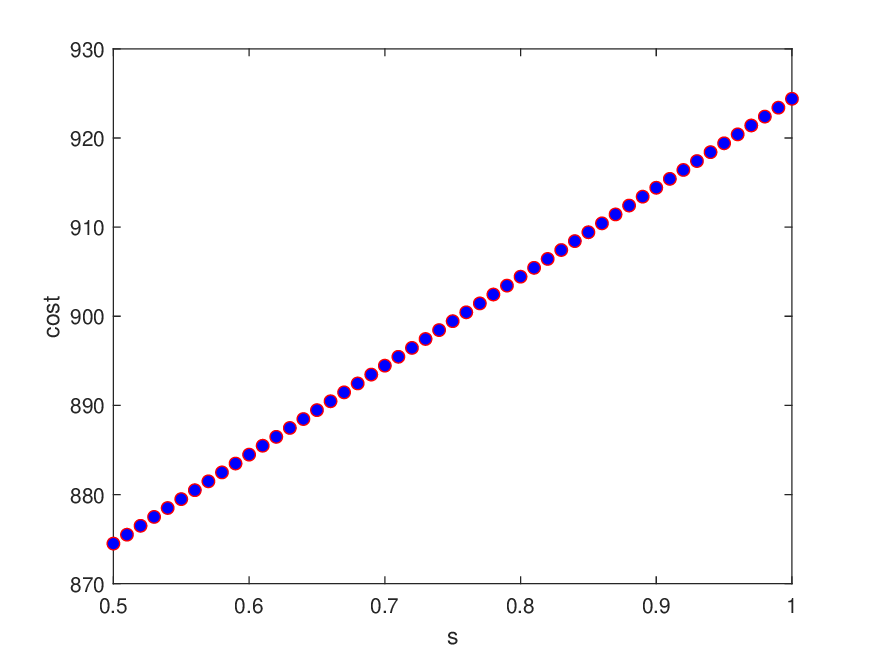} 
\vspace*{-1.0cm} 
\caption{{Optimal cost with varying $s$}}
\end{minipage}
\label{figsCost}
\end{figure}

\noindent We compared the optimal cost $923.9$ from our model by taking the similar parameters from Saadany \cite{Saadany2009} and compared it with their result their optimal cost $= 1635$ with $m=2$ and $n=1$. As mentioned earlier, the common practice is to consider combinations of $m$ and $n$ to get the optimal solution. Saadany \cite{Saadany2009} used various combinations of $m$ and $n$ to optimize cost and did not consider the combination $m=4$ and $n=3$, which turns out to be the optimal solution by our model. 
Also, Figure~3 shows that the inventory cost increases as $s$ increases with other parameters fixed, meaning we are sending more items to the repair depot, which results in increased cost. 

\par
\noindent{\bf Case II (Model 2):} When $s$ is not fixed and $R_1+R_2\geq D_r$\\
In this case, we control $s$ to send the remaining items from the repair depot to fulfil the demand. Table~\ref{table:exp-1a} presents the relationship among costs, new item cycle $n$, repair cycle $m$ and time cycle $T$ and $s$. 
 \begin{example}\label{ex3}\
In this example, we compare our model~2 \eqref{eqmathemodel2} with the model introduced in \cite[Model~I in Section~6.3]{Saadany2009}. This test example is introduced in \cite{Teunter2004}, and further analysis is conducted in \cite{Saadany2009} using their model~I. The following parameter settings are taken in the example, where $D_p=D_r=1000$, $p=r=0.8$, $q=1$,  $h_p=2$, $h_r=2$,  $S_p=20$ and $S_r=5$. To facilitate the comparison, we slightly adjust the parameters, setting $C_{1s}=C_{1r}=0$ (cost per unit of lost demand for new and repaired items). 
 \end{example}
The optimum solution occurs with the combination $m=1$, $n=7$ and $s=0.2$. The optimum cost is $924.07$ when $Q_p=143$ and $Q_r=1001$. We solved both models and obtained the optimum solution. After considering the parameters that provide us with the optimal solution, we created surfaces for both models, which are demonstrated in Figure \ref{figex43}. In this case, model~2 \eqref{eqmathemodel2} successfully approximates a better solution than the Model~I introduced in \cite[Model~I in Section~6.3]{Saadany2009}.

\begin{figure}[H]
\begin{center}
\hspace*{-2.5cm}
\includegraphics[width=80mm]{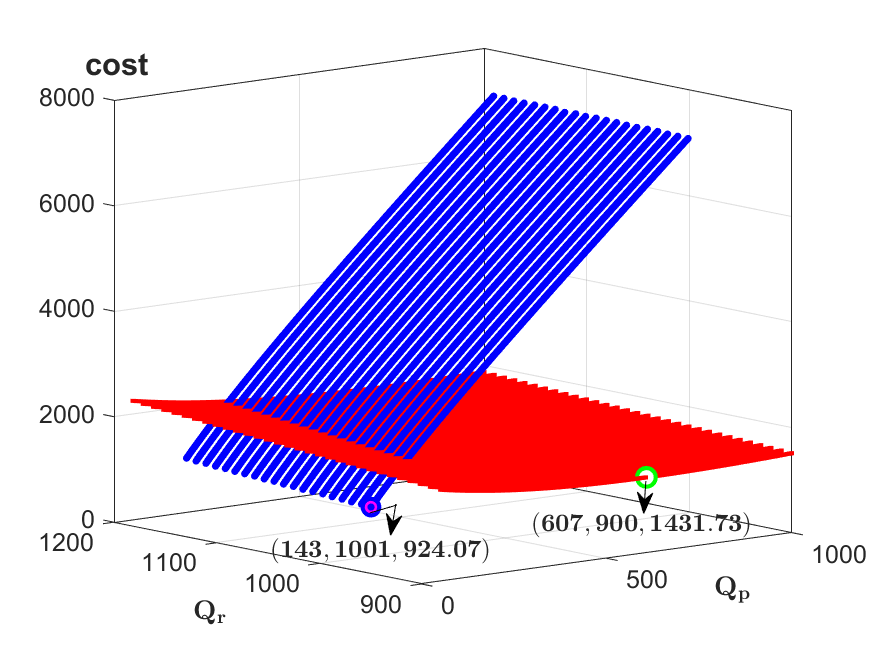}
\end{center}
\vspace*{-.8cm}
\caption{{Comparison of optimal solutions for Model \eqref{eqmathemodel2} and \cite[Model~I in Section~6.3]{Saadany2009}.}}
\label{figex43}
\end{figure}

\begin{example}\label{Swati}
This test example is taken from Swati et al \cite{Swati2021}. We assume that the demands in the primary and secondary markets are $D_p = 1500$ and $D_r = 2500$, respectively, over cycle $T$. The collection percentage of available returns from new items is $p = 0.8$, and after sorting, the recovery rate is $q = 0.8$. Similarly, the collection percentage of available returns from the secondary market is $r = 0.8$. After sorting, the recovery rate is $0 \leq s \leq 1$, where $s$ is a variable that depends on the number of items needed to fulfil the demand in the secondary market.

\noindent Assume that the holding costs per unit time for supply and repair depots are $h_p = 5$ and $h_r = 2$, respectively. The setup costs for new and repaired items are $S_p = 2400$ and $S_r = 1400$, respectively. We aim to determine the batch sizes $Q_p$ and $Q_r$, and the cycle numbers for repaired items $m$ and new items $n$, such that the average holding cost per unit time over cycle $T$ is minimised. Therefore, we seek to minimise the cost of Model 2 with constraints given by equation~(\ref{eqmathemodel2}) with the above parameter settings.
\end{example}

\noindent For the input settings in Example \ref{Swati}, the following results are obtained: total cost value $6372.5$ where $Q_p=1500$, $Q_r=1250$, $m=2$, $n=1$ and $s=0.77$. The results indicate that to minimise the cost to $6372.5$, the supply depot requires one cycle with a batch size of 1500 for new items and two cycles with 1250 for repaired items. The value $s=0.77$ indicates that after sorting, the recovery rate must be at least $77\%$ to meet the secondary market demand $D_r$. The algorithm's implementation takes approximately $1$ second of computational time. The total cost function with the optimum solution is demonstrated in Figure \ref{Ex44figsurf}. 

\begin{figure}[H]
\hspace{-1cm}
\begin{center}
\includegraphics[width=100mm]{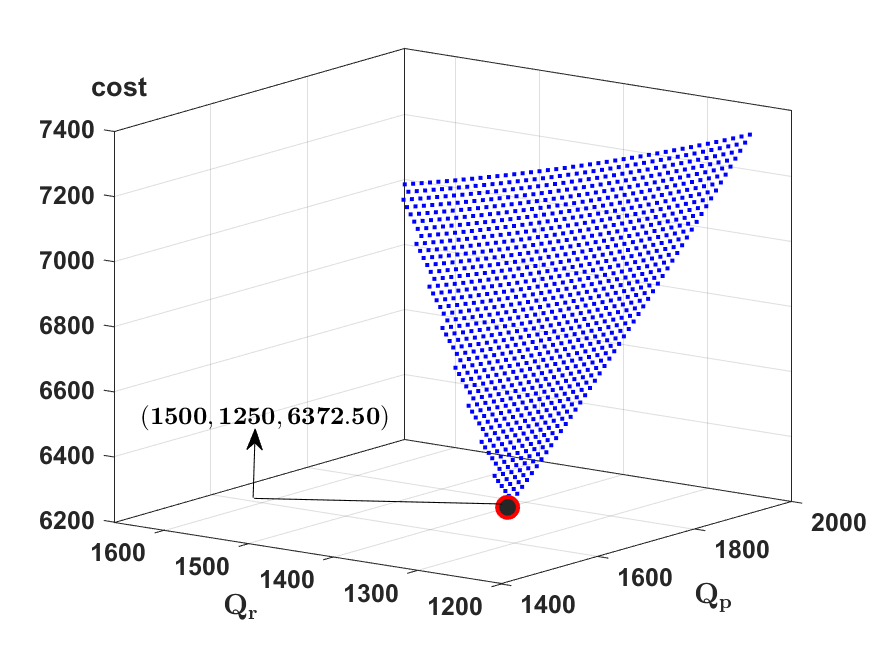} \\
\end{center}\hspace{-1cm}
\vspace{-1cm}
\caption{{ The total cost function with respect to $Q_p$ and $Q_r$ for Example \ref{Swati} when $m=2$, $n=1$, and $s=0.77$. 
}}
\label{Ex44figsurf}
\end{figure}

\subsection{Sensitivity Analysis} \label{sec4b}
A sensitivity analysis is conducted to examine the impact of parameter changes on the optimal solution. The results of this analysis are illustrated in Figures \ref{figSpSrCost}-\ref{fighrCost}. The observations are as follows.

\begin{enumerate}
\item [(i)] Figure~\ref{figSpSrCost}(a) demonstrates that keeping $S_r$ fixed results in no changes to the repair cycle $m$. However, when $S_p$ decreases, there is a significant increase in the number of procurement cycles $n$, which reduces inventory costs by $44\%$. Therefore, the optimal strategy is to increase the number of new procurement cycles as the setup cost $S_p$ decreases. Next, when $S_p$ is fixed and $S_r$ decreases, the number of repair cycles $m$ significantly increases, which reduces total inventory costs by $33\%$. Hence, the best strategy is to increase the number of repair cycles as the setup cost $S_r$ decreases, which is demonstrated in  Figure~\ref{figSpSrCost}(b).
\item [(ii)] Figure~\ref{fighpCost}(a) illustrates the impact on inventory cost when all parameters are fixed and the supply depot holding costs (\(h_p\)) is increased. When \(h_p = 5\), the optimal cost is \(6372.5\), achieved with the combination \(m = 2\) and \(n = 1\). If \(h_p\) increases to \(15\), the cost rises by $72.36\%$, with the new optimal combination being \(m = 3\) and \(n = 2\). The supply depot deals with both new and repaired items, so any changes to $h_p$ could potentially affect $m$ and $n$.
The cost continues to increase as \(h_p\) increases further, resulting in more cycles in the repair and supply depots. In this scenario, optimality is achieved with fewer cycles in repair and supply as \(h_p\) increases. Figure~\ref{fighpCost}(b) illustrates the impact of the other two variables, \(Q_p\) and \(Q_r\), and their relationship with the total cost as \(h_p\) changes.  
\item [(iii)] A similar trend is observed when the repair shop holding cost (\(h_r\)) increases. We increased \(h_r\) by up to $45\%$, demonstrating the effect in Figures~\ref{fighrCost}$(a)$ and $(b)$. Just as with \(h_p\), it is recommended to reduce the number of cycles in the repair depot as \(h_r\) increases while keeping the repair setup costs fixed.
\end{enumerate}

\begin{figure}[H]
\hspace{-1cm}
\begin{minipage}{100mm}
\begin{center}
\hspace*{-2.5cm}
\includegraphics[width=80mm]{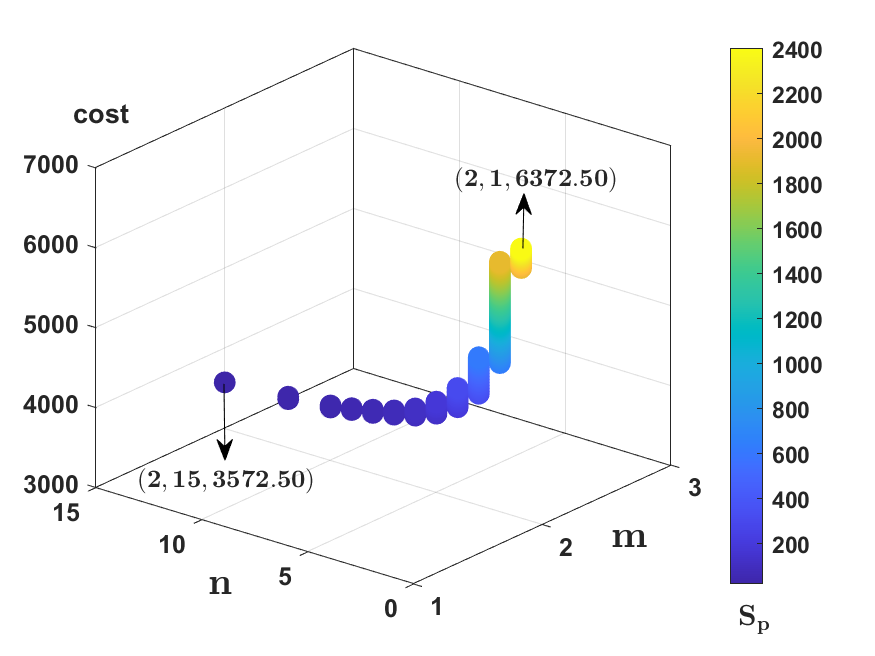} \\
\hspace*{-2.9cm}
{\scriptsize (a) Cost function with respect to $m$ and $n$\\ when $S_p$ increases.}
\end{center}
\end{minipage}
\hspace*{-3.5cm}
\begin{minipage}{100mm}
\begin{center}
\hspace*{0cm}
\includegraphics[width=80mm]{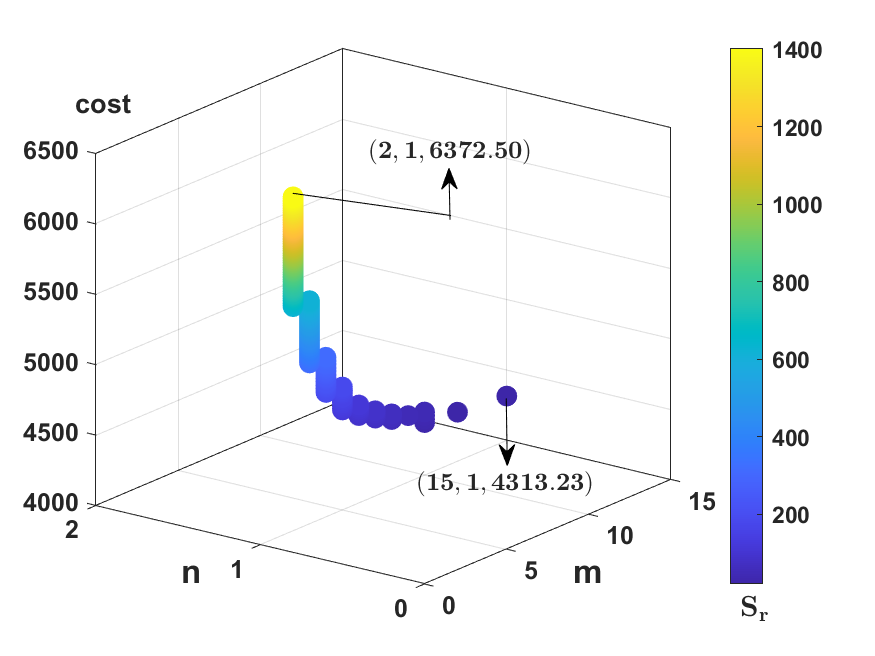} \\
\hspace*{-1.0cm}
{\scriptsize (b) Cost function with respect to $Q_p$ and $Q_r$\\ when $S_r$ increases.}
\end{center}
\end{minipage}
\caption{{ Effect of $(S_p)$ and $(S_r)$ on logistics inventory cost, repair and procurement cycles and batch sizes.
}}
\label{figSpSrCost}
\end{figure}
\begin{figure}[H]
\hspace{-1cm}
\begin{minipage}{100mm}
\begin{center}
\hspace*{-2.5cm}
\includegraphics[width=80mm]{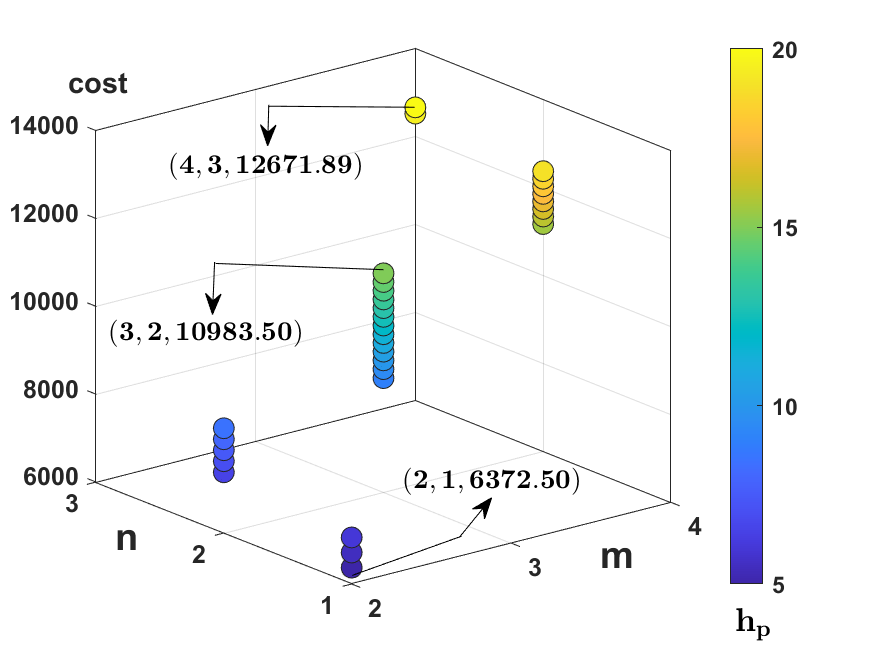} \\
\hspace*{-2.9cm}
{\scriptsize (a) Cost function with respect to $m$ and $n$\\ when $h_p$ increases.}
\end{center}
\end{minipage}
\hspace*{-3.5cm}
\begin{minipage}{100mm}
\begin{center}
\hspace*{0cm}
\includegraphics[width=80mm]{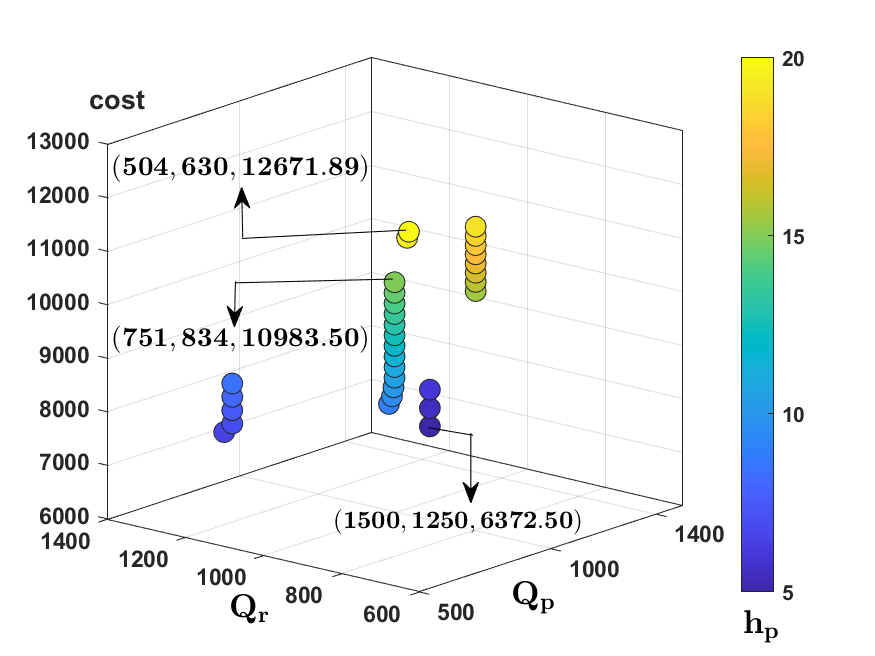} \\
\hspace*{-1.0cm}
{\scriptsize (b) Cost function with respect to $Q_p$ and $Q_r$\\ when $h_p$ increases.}
\end{center}
\end{minipage}
\caption{{Holding cost of supply depot $(h_p)$ effects logistics inventory cost, repair and procurement cycles and batch sizes.
}}
\label{fighpCost}
\end{figure}

\begin{figure}[H]
\hspace{-1cm}
\begin{minipage}{100mm}
\begin{center}
\hspace*{-2.5cm}
\includegraphics[width=80mm]{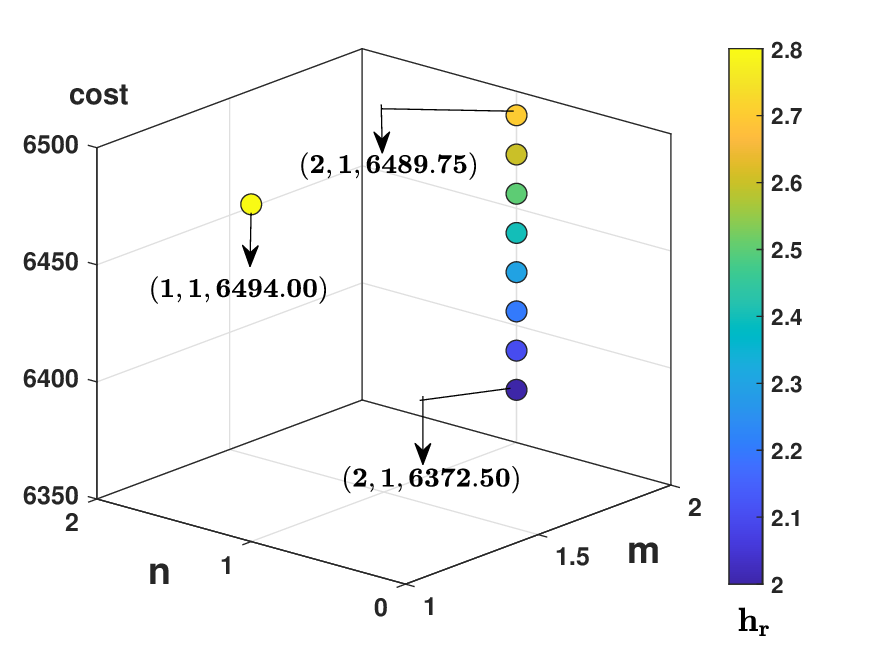} \\
\hspace*{-2.9cm}
{\small (a) Cost function with respect to $m$ and $n$\\ when $h_r$ increases.}
\end{center}
\end{minipage}
\hspace*{-3.5cm}
\begin{minipage}{100mm}
\begin{center}
\hspace*{0cm}
\includegraphics[width=80mm]{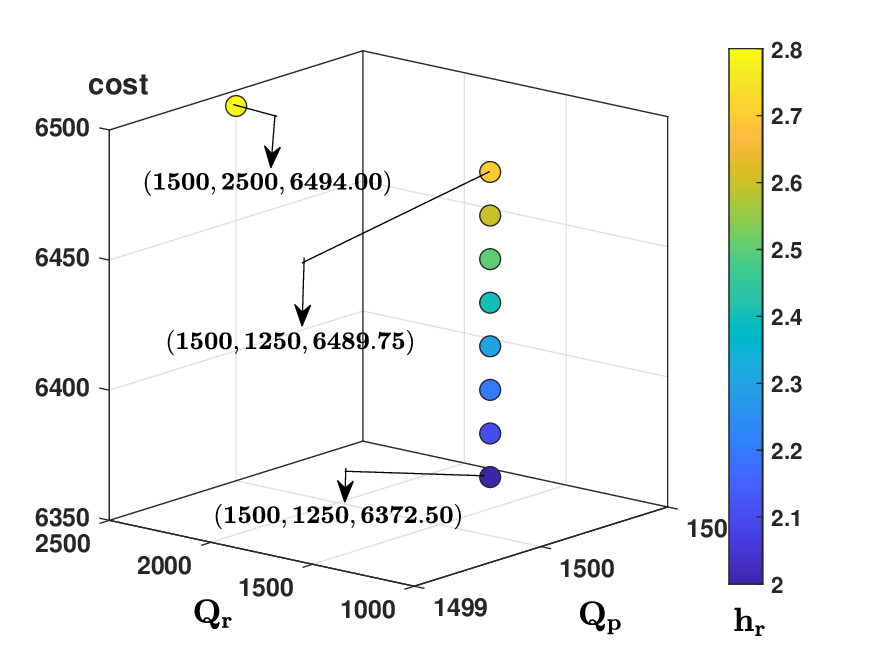} \\
\hspace*{-1.0cm}
{\scriptsize (b) Cost function with respect to $Q_p$ and $Q_r$\\ when $h_r$ increases.}
\end{center}
\end{minipage}

\caption{{Holding cost of repair depot $(h_r)$ effects logistics inventory cost, repair and procurement cycles and batch sizes. 
}}
\label{fighrCost}
\end{figure}

\newpage 
\subsection{Multiobjective Optimization Problems} \label{MOP}

\noindent Multiobjective optimization problems (MOPs) involve optimizing two or more conflicting objectives simultaneously, requiring a trade-off to find a set of optimal solutions known as the Pareto front. Scalarization techniques are a common approach to solving MOPs by transforming them into single-objective optimization problems. This transformation is achieved by combining the multiple objectives into a single scalar objective function, often using weights, which optimizes one objective while treating the others as constraints with specified bounds, and the achievement scalarizing function, which measures the distance from a solution to a reference point like the utopia point. These techniques facilitate the use of traditional optimization algorithms and provide decision-makers with a means to explore different trade-offs and prioritize conflicting objectives. In our analysis, we used two well-known scalarization techniques: the feasible-value-constraint approach \cite{BurKayRiz2017} and the Pascoletti Serafini approach \cite{PascolettiSerafini1984}. By applying both methods, we can evaluate the efficiency and robustness of our solutions, ensuring a complete approximation of the Pareto front. The descriptions of the approaches are given below.

\noindent
\textbf{\em The feasible-value-constraint approach \em\,:} For $\hat{x} \in X$, $X$ be a feasible set,
and
\[\hspace*{20mm} W := \left
\{ w \in \mathbb{R}^{\ell} \mid w_i > 0,\sum_{i=1}^{\ell} w_i=1
\right \}.\]
Then $w \in W$ and
\begin{equation}\label{property}
\hspace*{20mm} w_kf_k(\hat{x})=w_if_i(\hat{x}), \; \forall i=1,2,..., \ell, \;\; i \neq k.
\end{equation}
We assume a weight verifies \eqref{property}. We define the scalar problem as
\begin{equation} \label{nchep}
\mbox{($P_{\hat{x}}^k$)}\ \left\{\begin{array}{rl} \ds\min_{x\in  X} & \
\ w_{k}f_{k}(x),
\\[4mm]
\mbox{subject to} & \ \ w_if_i(x) \leq w_kf_k(\hat{x}), \;\;\mbox{i=1,2,...,$\ell$}, \;\; i \neq k.
\end{array}
\right.
\end{equation}

\noindent
\textbf{\em The Pascoletti-Serafini approach\em\,:} Another form of
the parameter-based scalarization approach introduced in \cite{PascolettiSerafini1984}, is also known as
\emph{goal-attainment method} (see, \cite{KayaJoydeep2011, Miettinen1999}).  For a given $w \in W$, the scalarization is
stated as follows.
\begin{equation} \label{P-S}
\hspace*{-5mm} \mbox{(PS)} \hspace{20mm}\;\; \ds\min_{(\alpha, x) \in \mathbb{R} \times X} \alpha, \;\;\;
\mbox{s.t.} \;\; w_i\,(f_i(x)-u_i) \leq  \alpha\,,\ \
\forall \;\;i=1,\ldots,\ell\,,
\end{equation}
where $\alpha \in \mathbb{R}$ is a new variable and $u$ is a utopia vector associated with Problem~($P$). 

\begin{example}\label{multexam}
	In this example, our objective is to solve multiobjective model 3 given by equation \eqref{eqmathemodel3} for the following parameter configuration: \(D_p=200\), \(D_r=50\), \(p=0.3\), \(q=0.9\), \(r=0.8\), \(S_p=4\), \(S_r=1\),  \(h_p=1\), \(h_r=1\), and \(c_w=1\). Our goal is to approximate the Pareto front of multiobjective model 3 using the scalarization methods $P_{\hat{x}}^k$ and  $PS$.
\end{example}

\noindent The point $(-70, -1000)$ is taken as the reference point (utopia vector). Equal grids $w_i$ are supplied in Algorithms~$1$ and $2$ (see, Appendix~\ref{appen}) to obtain the Pareto points, which are shown in Figure~\ref{multiexample}. The information about the number of grid points, subproblems, and solutions is presented in Table~\ref{tablemop}. The scalarization approach $P_{\hat{x}}^k$ approximated $278$ Pareto points, and one dominated point, and this dominated point is removed from the front. Due to the local solver trap within the neighbourhood, results approximate non-Pareto point. One can use the global solver to avoid dominated points; however, a global solver takes a long computation time compared to the local solver.

\begin{table}[ht]
\small
\caption{Numerical performance of $P_{\hat{x}}^k$ and $PS$.}  
\centering 
\begin{tabular}{||c| c| c| c| c| c||} 
\hline \hline 
Methods &   Number of  & Number of   & Time          & Pareto    & Non-Pareto  \\  
       &   grid points     & subproblems & in seconds    & points    &   points \\ [0.5ex] 
\hline 
$P_{\hat{x}}^k$ & 200       & 600       & 600       & 278        & 1     \\ [.5ex]\hline 
$PS$            & 200       & 200        & 52       & 200        & Nil     \\ 
\hline
\hline 
\end{tabular}
\label{tablemop} 
\end{table}

 \noindent \textbf{Analysis:} In Figures~\ref{multiexample}(a) and \ref{multiexample}(b), the objective function $f_1$ represents total inventory costs, with its range under the current parameter settings being between $\$48.50$ and $\$69.60$. On the other hand, $f_2$ represents waste disposal costs, ranging from $\$6.00$ to $\$46.00$. It's evident that $f_1$ has the lowest cost at $\$48.50$, while waste cost is at its highest at $\$46.00$. If one aims to decrease inventory costs, there's a need to compromise on waste costs. For the interval of inventory costs from $\$59.00$ to $\$69.60$, waste costs remain constant at $\$6.00$. This occurs when $s=1$ within this interval. The Pareto front is a continuous front, with some gaps due to algorithmic challenges in approximating certain points on the front. We utilise the Pascoletti-Serafini approach ($PS$) and the objective constraint approach $P_{\hat{x}}^k$, two well-known scalarisation approaches, to solve the multiobjective optimisation problem. General algorithms for both methods are detailed in the appendix section (see Appendix~\ref{appen}). While the algorithm employs multiple solvers, solving the problem proves challenging due to its mixed-integer multiobjective nature \cite{BurKayRiz2019}. We found that the local solver 'Bonmin' performs well in approximating Pareto solutions.
    
\noindent Waste and inventory costs largely vary depending on $s$, as illustrated in Figures \ref{multiexample}(c) and \ref{multiexample}(d). When $s$ is small, waste costs increase, indicating that not all items are considered for repair even if they meet the conditions for repair. However, waste costs decrease when $s=1$, meaning all repairable items are considered for repair.

\begin{figure}[H]
\hspace{-1.1cm}
\begin{minipage}{80mm}
\begin{center}
\hspace*{0cm}
\includegraphics[width=80mm]{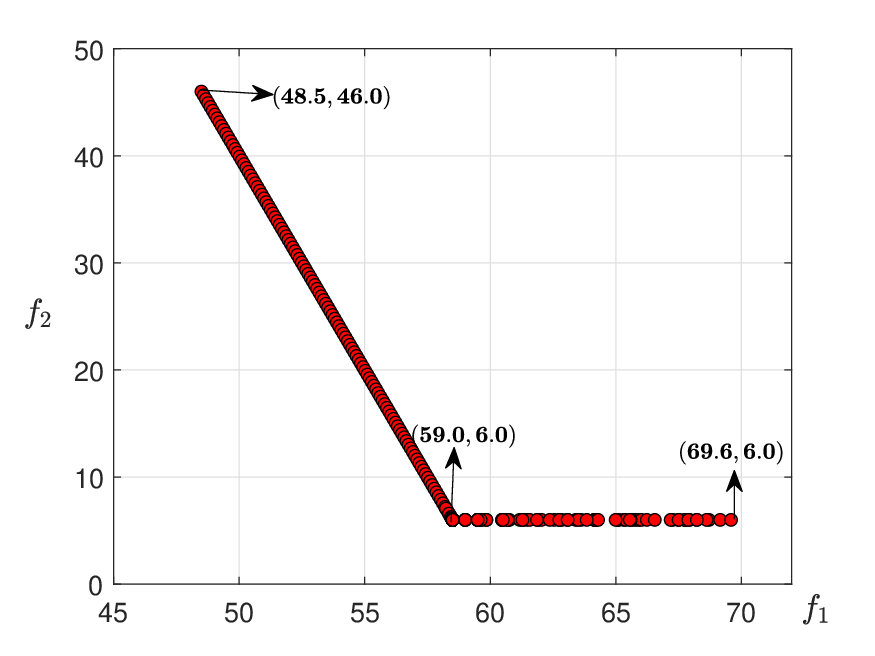} \\
{\scriptsize (a) Pareto points found with $P_{\hat{x}}^k$.}
\end{center}
\end{minipage}
\hspace{-0.7cm}
\begin{minipage}{80mm}
\begin{center}
\hspace*{0cm}
\includegraphics[width=80mm]{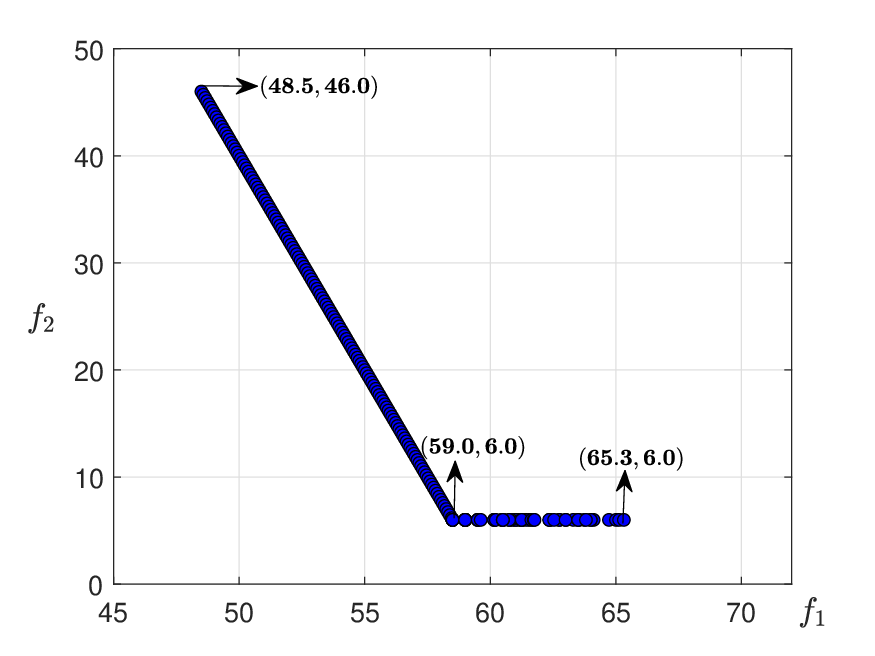} \\
{\scriptsize (b) Pareto points found with $PS$.}
\end{center}
\end{minipage}
\\[2mm]
\hspace*{-1.1cm}
\begin{minipage}{80mm}
\begin{center}
\hspace*{0cm}
\includegraphics[width=80mm]{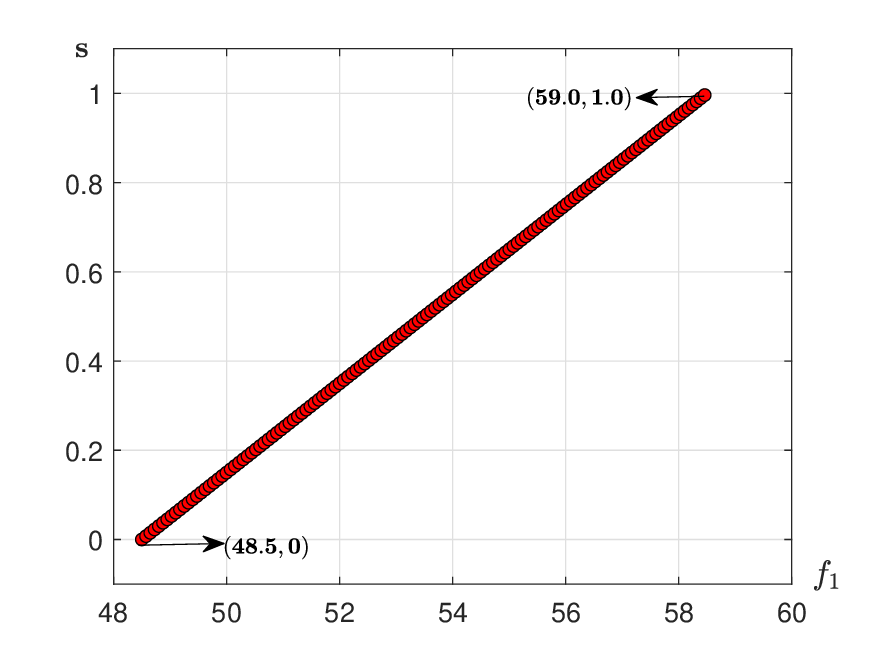} \\
{\scriptsize (c) Variation of $f_1$ for $0\leq s \leq 1$.}
\end{center}
\end{minipage}
\hspace{-0.6cm}
\begin{minipage}{80mm}
\begin{center}
\hspace*{0cm}
\includegraphics[width=80mm]{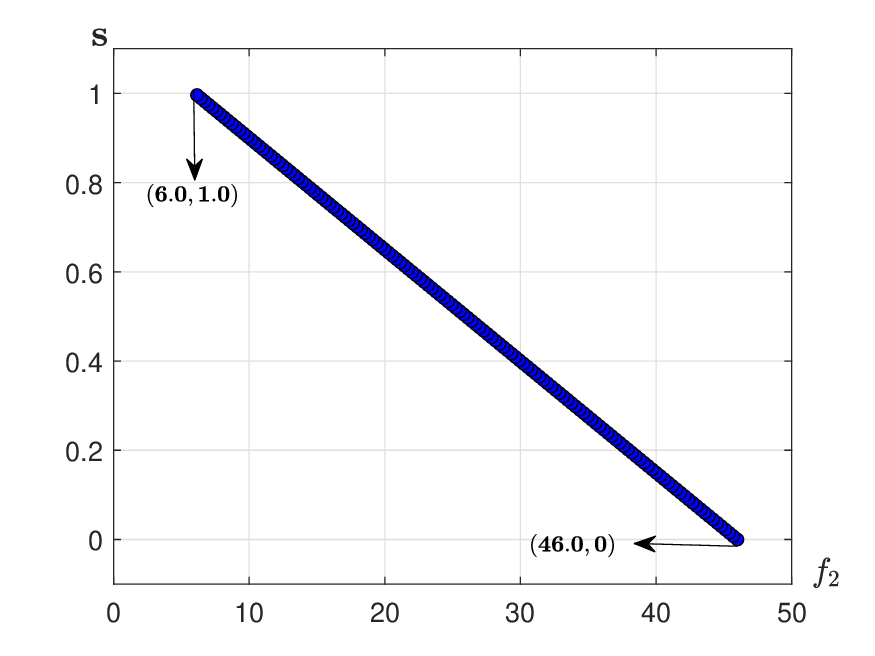}\\
{\scriptsize (d) Variation of $f_2$ for $0\leq s \leq 1$.}
\end{center}
\end{minipage}

\caption{{Equal grids are provided in Algorithm.
            (Red) circles indicate the Pareto front obtained by $P_{\hat{x}}^k$ in (a), (blue) circles depict Pareto points found by ($PS$) in (b), variation of $f_1$ and $f_2$ over $0 \leq s \leq 1$ is shown in (c) and (d).}}
\label{multiexample}
\end{figure}
\newpage

\section{Conclusion}\label{Con}
This study explores a reverse logistics inventory model dealing with different primary and secondary market demands. While existing literature has examined separate markets for new and repaired items, none have addressed the issue of incorporating different types of parameters (integers and continuous) together. In our paper, we developed three models with distinct constraints and used Mixed Integer Programming to find the optimal solution. Model~3 incorporated environmental factor disposal cost as a secondary objective function. Additionally, we utilized two scalarization techniques to identify the Pareto front of the multiobjective problem. This gives manufacturers the flexibility to create their own optimal policy by balancing holding costs and repair costs. 

\noindent For future work, we suggest expanding on the deterministic and stochastic Economic Production Quantity (EPQ) model for reverse logistics by using MIP and environmental factors such as GHG and energy consumption. 

\noindent \textbf{\large {Acknowledgment:}} 
\noindent We acknowledge Chatgpt(3.5) AI tools, which are used to rephrase, edit and polish the author's written text for spelling, grammar, or general style. 

\section{\small Funding and/or Conflicts of interests/Competing interests}
 The authors declare that no funds, grants, or other support were received during the preparation of this manuscript. The authors have no relevant financial or non-financial interests to disclose.



\section{Appendix}\label{appen}

\subsection{$C(Q_p, Q_r, m, n)$ formulation.}
\noindent The calculations for the cost function \eqref{eqcost1a} based on Figure~\ref{figinventoryModel1} are as below.
\begin{eqnarray*} \label{eqcost}
Area~\mbox{KLM} &=& \frac{\mbox{KL.LM}}{2} = \frac{Q_p^2}{2D_p}.\\
Area~\mbox{NOS} &=& \frac{\mbox{NO.OS}}{2} = \frac{Q_r^2}{2D_r} .\\
Area~\mbox{XYF} &=& \frac{\mbox{XF.YF}}{2} = \frac{T_p^2pqD_p}{2}=\frac{R_1T_p^2}{2}.\\
Area~\mbox{ABC} &=& \frac{\mbox{AC.BC}}{2} = \frac{Q_r^2 srD_r}{2D_r^2} = \frac{Q_r^2 R_2}{2D_r^2}. \\
Area~\mbox{BCDE} &=& {\mbox{BC}}\times{\mbox{BD}} = \frac{Q_r^2}{D_r}\bigg(1-\frac{R_2}{D_r}\bigg).\\
Area~\mbox{DEGF} &=& {\mbox{DE}}\times\mbox{{DF}} =\frac{Q_r}{D_r}(YF-YD) = \frac{Q_r}{D_r}\bigg(T_pR_1-(Q_r+(m-1)Q_r\bigg(1-\frac{R_2}{D_r}\bigg)\bigg).
\end{eqnarray*}
Adding all areas, the inventory cost function is given as 
\begin{multline*}
C(Q_p, Q_r, m, n) = \frac{1}{T}\Bigg[ mS_n+nS_p + \left(\frac{Q_pT_p}{2}+\frac{Q_rT_r}{2}\right)h_p +\Bigg[ \frac{R_1T_p^2}{2} + \frac{R_2T_rQ_r}{2D_r} \\ +\frac{(m-1)}{2}T_r Q_r\left(1-\frac{R_2}{D_r}\right) 
+T_r\left(T_pR_1-Q_r-(m-1)\left(1-\frac{R_2}{D_r}\right)Q_r\right)\Bigg] h_r\Bigg],
\end{multline*}
where, $T=T_p+T_r$, $R_1 = pqD_p$ and $R_2 = srD_r$.
\newpage 
\subsection{Algorithms}
\noindent We present Algorithms $1$ and $2$ which are used in Example-\ref{multexam} to generate the Pareto fronts for the scalarization approaches ($P_{\hat{x}}^k$) and $PS$, respectively. 

\noindent \textbf{Algorithm 1} 
\begin{description}
	\item[Step $\mathbf{1}$] { \textbf{(Input)}} \\ Set all parameters that stated in Example-\ref{multexam}.
	\item[Step  $\mathbf{2}$] {\textbf{(Determine the individual minima)}}\\
	Solve  $\ds \min f_i, \;\; i=1,2$, subject to the constraints given in model \eqref{eqmathemodel3} for parameter settings stated in Example-\ref{multexam}, and obtain solutions $\left(\bar{Q}_{p_i}, \bar{Q}_{r_i}, \bar{m}_i, \bar{n}_i, \bar{s}_i\right)$, for $i=1,2$, respectively.	
 
	\item[Step $\mathbf{3}$] { \textbf{(Generate weighted parameters)}} \\
	In this step, we generate weight parameters $w_i$. Following the weight generation technique given in \cite[Step 3 of Algorithm 3]{BurKayRiz2019}, we generate weights considering the boundary solution obtained in Step~1. 
 \item[Step $\mathbf{4}$]{ \textbf{(Locate $\hat{x}$)}} \\
 Find $\hat{x}$ that solves an additional problem $\ds\min_{x\in  X}\;\;\ w_{k}f_{k}(x),$ s.t. $w_if_i(x) = w_kf_k(\hat{x})$, for some $k$.
	\item[Step $\mathbf{5}$] Choose $w=(w_1,1-w_1)$, which generated from Step 3.
	\begin{description}
		\item[(a)] Find $\hat{x}_k=\left(Q_{p_k}, Q_{r_k}, m_k, n_k, s_k\right)$ that solves Problems ($P_{\hat{x}}^k$), $k=1,2.$
		\item[(b)] Determine weak efficient points : 
				\begin{description}
			\item[(i)] If $\bar{x}_1=\bar{x}_2$, then set $\bar{x}=\bar{x}_1$ (an efficient point)\\
			and store the efficient point.
			\item [(ii)] If $\bar{x}_1=\bar{x}_2$  does not hold, then any dominated point is discarded by comparing these $2$ solutions. \\
			Store non-dominated Pareto points.	
		\end{description}
	\end{description}
	\item[Step $\mathbf{6}$] (Output)\\
	All recorded points are Pareto points of Example-\ref{multexam}.\\
\end{description}
\newpage
\noindent \textbf{Algorithm 2}

\begin{description}
	\item[Steps $\mathbf{1}$-$\mathbf{3}$] Same as Algorithm 1.
 	\item[Step $\mathbf{4}$] Choose $w=(w_1,1-w_1)$ which generated from Step 3. \\ \hspace{10mm} Set utopia point $u=(-70, -1000)$.
	\begin{description}
		\item[(a)] Find $\hat{x}_k=\left(Q_{p_k}, Q_{r_k}, m_k, n_k, s_k\right)$ that solves Problems ($PS$).
		\item[(b)] Determine weak efficient points and record the points. 
	\end{description}
	\item[Step $\mathbf{5}$] (Output)\\
	All recorded points are Pareto points of Example-\ref{multexam}.\\
\end{description}

\subsection{Alternative flow model}
This is the alternative inventory flow diagram for new and returned (new and repaired) items. In this model, the flow of new and used items in the supply depot is represented in parallel as opposed to side by side, which is common in literature. The inventory cost can also be minimised by the cost function represented by equation (\ref{eqcost_parallel}). 
\begin{figure}[H]
\hspace{-1cm}
\begin{center}
\includegraphics[width=120mm]{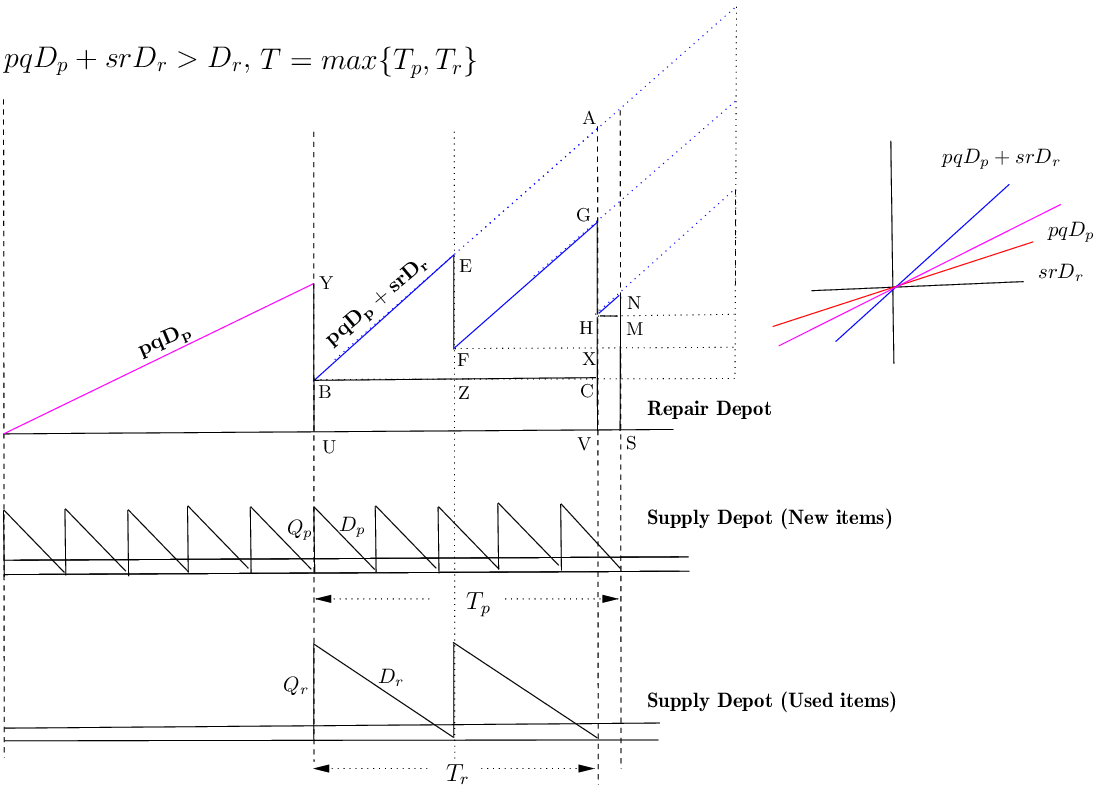} \\
\end{center}\hspace{-1cm}
\vspace{-1cm}
\caption{{The behaviour of inventory for produced, collected, used and repaired items over the interval $T$ in parallel flow.
}}
\label{Ex2figsurf}
\end{figure}

\noindent Mathematical expression of the inventory cost function is as follows.

\begin{multline} \label{eqcost_parallel}
f_A(Q_p, Q_r, m, n, s) = \frac{1}{T}\Bigg[ mS_n+nS_p + \left(\frac{Q_pT_p}{2}+\frac{Q_rT_r}{2}\right)h_p + \Bigg[\frac{1}{2m}{T_r^2}(R_1+R_2) \\ + \frac{(m-1)}{2}T_r Q_r\left(\frac{R_1+R_2}{D_r}-1\right)+\frac{1}{2}(T_p-T_r)^2(R_1+R_2)+(T_p-T_r)\Big(T_r(R_1+R_2)\\+R_1T_p-(m+1)Q_r\Big)+(R_1T_p-Q_r)T_r\Bigg]h_r\Bigg],
\end{multline}
where, $T=T_p+T_r$, $R_1 = pqD_p$ and $R_2 = srD_r$.

\noindent \textbf{Author Contributions}: 
Both authors contributed to the study conception and design. Material preparation, data collection and analysis were performed by Dr Indu Bala Wadhawan and Dr Mohammed Mustafa Rizvi. 
\end{document}